\documentclass[12pt]{article}

\begin{document}

\def\gint#1{{\bf [}#1{\bf ]}}
\def\A{{\rm{\bf A}}}
\def\modd#1{\vert#1\vert}
\def\Hs{\heartsuit}
\def\div{{\rm div}}
\def\N{{\bf N}}
\def\Arr{{\cal R}}
\def\nrm#1{\Vert#1\Vert}
\def\dq{\lq\lq}
\def\fstar{f^*}
\def\rmk{{\emph{Remark}}}
\def\quote#1{\lq\lq#1''}
\def\Rd{{{\bf R}^d}}
\def\Upperhalfspace{{{\bf R}^{d+1}_+}}
\def\Dee{{\cal D}}
\def\Bee{{\cal B}}
\def\Cinf{{\cal C}^\infty_0(\Rd)}
\def\Lwunloc{L^1_{loc}(\Rd)}
\def\Lwunlocmu{L^1_{loc}(\mu)}
\def\Lwunlocv{L^1_{loc}(v)}
\def\fn#1#2{\footnote{#1}{#2}}
\def\fnno#1#2{\footnote{$ ^{#1}$}{#2}}
\def\sqrterm{{\nrm{a_Q(f)}^2_2\over\modd Q}}
\def\Ainf{A_\infty}
\def\Gee{{\cal G}}
\def\Eff{{\cal F}}
\def\M{{\cal M}}
\def\W{{\cal W}}
\def\I{{\cal I}}
\def\Eee{{\cal E}}
\def\Sss{{\cal S}}
\def\Uball{\{x:\ \modd x\leq1\}}
\def\Ca{Calder\'on}
\def\Pf{{\bf Proof.}\ \ }
\def\AN{{\cal A}_N}
\def\Ckay{{\cal C}^k(\Rd)}
\def\Sf#1#2#3{\left(\int\limits_{\modd{x-t}<#3 y}\modd{#1 * #2_y(t)}^2\,
{dt\, dy\over y^{d+1}}\right)^{1/2}}
\def\Sq#1#2#3{S_{#2,#3}(#1)}
\def\R{{\bf R}}
\def\bigmodd#1{\left\vert#1\right\vert}
\def\bignrm#1{\left\Vert#1\right\Vert}
\def\stars{\smallskip*\hfil*\hfil*\smallskip}
\def\Bispace{{\R^{d_1}\times\R^{d_2}}}
\def\leaderfill{\leaders\hbox to 1em{\hss.\hss}\hfil}
\def\topline#1{\headline={\ifnum\pageno>1 wilson\hfil #1\else\hfil\fi}}
\def\bigbrace#1{{\left\lbrace#1\right\rbrace}}
\def\Calpha{{{\cal C}_\alpha}}
\def\CalphaM{{{\cal C}_{\alpha,M}}}
\def\CalphaO{{{\cal C}_{\alpha,0}}}
\def\CalphaMO{{{\cal C}_{\alpha,0,M}}}
\def\qed{{$\clubsuit$}}
\def\Compnum{{\bf C}}
\def\Irreg{{\cal IR}}
\def\rmk{\emph{Remark}}
\def\C{\Compnum}
\def\measupp{{\{x:\ f(x)\not=0\}}}
\newtheorem{theorem}{Theorem}
\newtheorem{lemma}{Lemma}
\newtheorem{definition}{Definition}
\newtheorem{corollary}{Corollary}

\title{Growth and decay of H\"older moduli}\author{Michael Wilson\\
Department of
Mathematics\\University of Vermont\\Burlington, Vermont\ \
05405\\jmwilson@uvm.edu}\date{}

\maketitle

\begin{abstract}
If $f:\Rd\to\C$ is bounded and $f$'s H\"older $\alpha$-modulus of continuity grows no faster than $(1+\modd x)^M$ ($M\geq0$) then, for every $\epsilon>0$, there is a $\beta>0$ such that $f$'s H\"older $\beta$-modulus grows no faster than $(1+\modd x)^{\epsilon}$. We use this easy fact to show that, if $\modd f$ decays as fast as $(1+\modd x)^{-R}$ (for $R>0$) and $f$'s $\alpha$-H\"older modulus grows no faster than $(1+\modd x)^M$, then, for every $0\leq R'< R$, there is a $\beta>0$ such that $f$'s $\beta$-H\"older modulus decays as fast as $(1+\modd x)^{-R'}$. We apply this to strengthen a result of Coifman and Meyer on almost-orthogonality of vaguelet families and to derive other useful facts about vaguelets and vaguelet-like functions.
\end{abstract}

In \cite[p.~56, Theorem 2]{WavCoifMeyer} Coifman and Meyer proved

\begin{theorem}\label{almostorthog}Let $0<\alpha<\epsilon$ and $C_1,C_2>0$. Let $\{f^{(Q)}\}_{Q\in\Dee_d}$ be a family of functions indexed over the dyadic cubes $\Dee_d$ in $\Rd$, where each $f^{(Q)}$ has integral 0 and: i) for all $x\in\Rd$,
\begin{equation}\label{decay}\modd{f^{(Q)}(x)}\leq C_1(1+\modd x)^{-d-\epsilon};\end{equation}
ii) for all $x,x'\in\Rd$, 
\begin{equation}\label{beta}\modd{f^{(Q)}(x)-f^{(Q)}(x')}\leq C_2\modd{x-x'}^\alpha.\end{equation}
For each $Q\in\Dee_d$, let $x_Q$ be $Q$'s center and $\ell(Q)$ its sidelength, and define
$$g^{(Q)}(x):={f^{(Q)}((x-x_Q)/\ell(Q))\over\modd Q^{1/2}},$$
where $\modd Q$ means $Q$'s Lebesgue measure. There is a constant $A$, depending only on $\epsilon$, $\alpha$, $C_1$, $C_2$, and $d$, so that, for all finite sequences $\{\lambda_Q\}_{Q\in\Dee_d}\subset\C$,
\begin{equation}\label{calledao}\int\bigmodd{\sum \lambda_Q g^{(Q)}}^2\,dx\leq A\sum\modd{\lambda_Q}^2.\end{equation}
\end{theorem}

\rmk. The hypothesis $\alpha<\epsilon$ is unnecessary (see Lemma \ref{independent} below).\medskip

Coifman and Meyer call the functions $g^{(Q)}$ \emph{vaguelets}. Theorem \ref{almostorthog} shows that $\{g^{(Q)}\}_{Q\in\Dee_d}$ is a \emph{Bessel family}; property (\ref{calledao}) is sometimes called \emph{almost-orthogonality}. Theorem \ref{almostorthog} is well known. Vaguelets are familiar objects in frame theory and signal processing (\cite{AbramSilv}\cite{CaiPlanning}\cite{DidJafPip}\cite{Jia}\cite{LiuYY}) if for no other reason than that they are what you get when you apply singular integral operators to wavelets.

In this paper we prove

\begin{theorem}\label{counterintuitive}Let $D_1,R,\alpha$ be positive numbers and $M\geq0$. For every $0\leq R'< R$ there are positive constants $D_2$ and $\beta$, depending only on $D_1, \alpha, R, R'$, and $M$, such that the following holds: If $f:\Rd\to\C$ satisfies
\begin{equation}\label{decayinf}\modd{f(x)}\leq (1+\modd x)^{-R}\end{equation}
for all $x\in\Rd$ and
\begin{equation}\label{growholder}\modd{f(x)-f(x')}\leq D_1\modd{x-x'}^\alpha(1+\max(\modd x,\modd{x'}))^M\end{equation}
for all $x,x'\in\Rd$, then
\begin{equation}\label{decayholder}\modd{f(x)-f(x')}\leq D_2\modd{x-x'}^\beta(1+\min(\modd x,\modd{x'}))^{-R'}\end{equation}
for all $x,x'\in\Rd$.
\end{theorem}

\rmk. The switch from \lq max' in (\ref{growholder}) to \lq min' in (\ref{decayholder}) is not a typo. The theorem says: If $f$ decays as fast as $(1+\modd x)^{-R}$ and its $\alpha$-H\"older modulus grows no faster than a polynomial then, for some $\beta>0$, $f$'s $\beta$-H\"older modulus decays almost as fast as $(1+\modd x)^{-R}$. The relevance of this to Theorem \ref{almostorthog} is as follows. Put $R:=d+\epsilon$ and suppose that every $f^{(Q)}$ in $\{f^{(Q)}\}_{Q\in\Dee_d}$ has integral 0 and satisfies (\ref{decayinf}) and (\ref{growholder}). By Theorem \ref{counterintuitive}, every $f^{(Q)}$ satisfies (\ref{decayinf}) for this $R$ and (\ref{decayholder}) for $R':=0$ and some $\beta>0$.  Theorem \ref{almostorthog} follows by our preceding remark. Theorem \ref{counterintuitive} also has this nice corollary.

\begin{corollary}\label{vaguelet}Let $\epsilon,C_1,C_2>0$ and $M\geq0$. Let $\{f^{(Q)}\}_{Q\in\Dee_d}$ be a family of differentiable functions indexed over the dyadic cubes $\Dee_d$ in $\Rd$, where each $f^{(Q)}$ has integral 0 and, for all $x\in\Rd$,
$$\modd{f^{(Q)}(x)}\leq C_1(1+\modd x)^{-d-\epsilon}$$
and
$$\modd{\nabla f^{(Q)}(x)}\leq C_2(1+\modd x)^M.$$
If the functions $g^{(Q)}$ are defined as in Theorem \ref{almostorthog} then (\ref{calledao}) holds for some $A$ depending only on $\epsilon,C_1,C_2,M$, and $d$.
\end{corollary}

Theorem \ref{counterintuitive} follows from a simple lemma which may be of independent interest. We believe it's well known, but we haven't seen it anywhere.

\begin{lemma}\label{independent}Let $A,\alpha,M\geq0$ and let $f:\Rd\to\C$ satisfy: a) for all $x\in\Rd$, $\modd{f(x)}\leq 1$; b) for all $x,x'\in\Rd$, 
$$\modd{f(x)-f(x')}\leq A\modd{x-x'}^\alpha(1+\max(\modd x,\modd{x'}))^M.$$
Then, for all $0<\rho\leq1$ and $x,x'\in\Rd$,
$$\modd{f(x)-f(x')}\leq \max(2,A)\modd{x-x'}^{\rho\alpha}(1+\max(\modd x,\modd{x'}))^{\rho M}.$$
\end{lemma}

{\bf Proof of Lemma \ref{independent}.} Trivial. If $\modd{x-x'}^\alpha(1+\max(\modd x,\modd{x'}))^M\geq1$ then
\begin{eqnarray*}\modd{f(x)-f(x')}&\leq& 2\\
&\leq&2\left(\modd{x-x'}^\alpha(1+\max(\modd x,\modd{x'}))^M\right)^\rho\\
&=&2\modd{x-x'}^{\rho\alpha}(1+\max(\modd x,\modd{x'}))^{\rho M}.\end{eqnarray*}
If $\modd{x-x'}^\alpha(1+\max(\modd x,\modd{x'}))^M\leq1$ then
\begin{eqnarray*}\modd{f(x)-f(x')}&\leq&A\modd{x-x'}^\alpha(1+\max(\modd x,\modd{x'}))^M\\
&\leq&A\left(\modd{x-x'}^\alpha(1+\max(\modd x,\modd{x'}))^M\right)^\rho\\
&=&A\modd{x-x'}^{\rho\alpha}(1+\max(\modd x,\modd{x'}))^{\rho M},\end{eqnarray*}
proving the lemma.\medskip

\rmk. Lemma \ref{independent} justifies the remark following the statement of Theorem \ref{almostorthog}.

{\bf Proof of Theorem \ref{counterintuitive}.}  By Lemma \ref{independent} we can assume that 
$$\modd{f(x)-f(x')}\leq D_1\modd{x-x'}^\alpha(1+\max(\modd x,\modd{x'}))^{R/2}$$
for all $x,x'\in\Rd$. We assume $\modd x\leq\modd{x'}$ and consider two cases:  i) $\modd{x-x'}^\alpha\geq (1+\modd x)^{-R}$ and ii) $\modd{x-x'}^\alpha\leq (1+\modd x)^{-R}$. For case i) we have
\begin{eqnarray*}\modd{f(x)-f(x')}&\leq&\modd{f(x)}+\modd{f(x')}\\
&\leq& (1+\modd x)^{-R}+ (1+\modd{x'})^{-R}\\
&\leq& 2(1+\modd x)^{-R}\\
&\leq&2(1+\modd x)^{-R/2}\\
&\leq&2\modd{x-x'}^{\alpha/2}.\end{eqnarray*}
For case ii) we have $\modd{x-x'}\leq1$, so that $(1+\modd{x'})\leq 2(1+\modd x)$, implying
\begin{eqnarray*}\modd{f(x)-f(x')}&\leq&D_1\modd{x-x'}^\alpha(1+\modd{x'})^{R/2}\\
&=&D_1\modd{x-x'}^{\alpha/2}\modd{x-x'}^{\alpha/2}2^{R/2}(1+\modd x)^{R/2}\\
&\leq&D_12^{R/2}\modd{x-x'}^{\alpha/2}.\end{eqnarray*}
Conclusion: there are positive constants $\rho$ and $D_3$ such that, for all $x,x'\in\Rd$,
$$\modd{f(x)-f(x')}\leq D_3\modd{x-x'}^\rho.$$

Define 
$$\gamma:={\rho R'\over R},$$
which satisfies $0\leq\gamma<\rho$.

Take $\modd x\leq\modd{x'}$ again and consider two cases: a) $\modd{x-x'}^\rho\geq(1+\modd x)^{-R}$ and b) $\modd{x-x'}^\rho\leq(1+\modd x)^{-R}$. For case a) we observe that
\begin{eqnarray*}\modd{x-x'}^{\rho-\gamma}&\geq&(1+\modd x)^{-R(1-\gamma/\rho)}\\
&=&(1+\modd x)^{-R}(1+\modd x)^{R\gamma/\rho}\\
&=&(1+\modd x)^{-R}(1+\modd x)^{R'},\end{eqnarray*}
implying
$$(1+\modd x)^{-R}\leq \modd{x-x'}^{\rho-\gamma}(1+\modd x)^{-R'}.$$
Since $(1+\modd{x'})^{-R}\leq(1+\modd x)^{-R}$, we also have $\modd{x-x'}^\rho\geq(1+\modd{x'})^{-R}$, and a similar argument yields
$$(1+\modd{x'})^{-R}\leq \modd{x-x'}^{\rho-\gamma}(1+\modd{x'})^{-R'}\leq  \modd{x-x'}^{\rho-\gamma}(1+\modd{x})^{-R'}.$$
Therefore, for case a),
\begin{eqnarray*}\modd{f(x)-f(x')}&\leq&\modd{f(x)}+\modd{f(x')}\\
&\leq&(1+\modd x)^{-R}+(1+\modd{x'})^{-R'}\\
&\leq&\modd{x-x'}^{\rho-\gamma}(1+\modd x)^{-R'}+\modd{x-x'}^{\rho-\gamma}(1+\modd x)^{-R'}\\
&=&2\modd{x-x'}^{\rho-\gamma}(1+\modd x)^{-R'}.\end{eqnarray*}

For case b) we have:
\begin{eqnarray*}\modd{f(x)-f(x')}&\leq& D_3\modd{x-x'}^\rho\\
&=&D_3\modd{x-x'}^{\rho-\gamma}\modd{x-x'}^\gamma\\
&\leq&D_3\modd{x-x'}^{\rho-\gamma}(1+\modd x)^{-R\gamma/\rho}\\
&=&D_3\modd{x-x'}^{\rho-\gamma}(1+\modd x)^{-R'},\end{eqnarray*}
giving us (\ref{decayholder}) with $\beta:=\rho-\gamma$ and $D_2:=\max(2,D_3)$. Theorem \ref{counterintuitive} is proved. \medskip

In wavelet/frame theory, (\ref{decayinf}) usually shows up with $R=d+\epsilon$. That makes it natural to consider four possible sets of hypotheses on a function $f:\Rd\to\C$.\medskip

$(I)$. There are positive constants $\epsilon$, $\alpha$, and a non-negative $M$, such that:\smallskip

For all $x\in\Rd$, 
$$\modd{f(x)}\leq (1+\modd x)^{-d-\epsilon}$$
and, for all $x,x'\in\Rd$, 
$$\modd{f(x)-f(x')}\leq \modd{x-x'}^\alpha\left(1+\max(\modd x,\modd{x'})\right)^M.$$\smallskip

$(II)$ There are positive constants $\epsilon$ and $\alpha$ such that:\smallskip

For all $x\in\Rd$, 
$$\modd{f(x)}\leq (1+\modd x)^{-d-\epsilon}$$
and, for all $x,x'\in\Rd$, 
$$\modd{f(x)-f(x')}\leq \modd{x-x'}^\alpha.$$\smallskip

$(III)$ There are positive constants $\epsilon$ and $\alpha$ such that:\smallskip

For all $x\in\Rd$, 
$$\modd{f(x)}\leq (1+\modd x)^{-d-\epsilon}$$
and, for all $x,x'\in\Rd$, 
$$\modd{f(x)-f(x')}\leq \modd{x-x'}^\alpha\left((1+\modd x)^{-d-\epsilon}+(1+\modd{x'})^{-d-\epsilon}\right).$$
\smallskip

$(IV)$ There are positive constants $\epsilon$ and $\alpha$ such that:\smallskip

For all $x\in\Rd$, 
$$\modd{f(x)}\leq (1+\modd x)^{-d-\epsilon}$$
and, for all $x,x'\in\Rd$, 
$$\modd{f(x)-f(x')}\leq \modd{x-x'}^\alpha\left((1+\modd x)^{-d-\epsilon-\alpha}+(1+\modd{x'})^{-d-\epsilon-\alpha}\right).$$
\smallskip

It's trivial that (neglecting unimportant multiplicative factors) $(IV)\Rightarrow(III)\Rightarrow(II)\Rightarrow(I)$. Theorem \ref{counterintuitive} shows that, with the same proviso, $(I)\Rightarrow(III)$. In \cite[Definition 6.2, Lemma 6.2]{Wilsonbook} it is shown that $(III)$ implies $(IV)$. Therefore $(I)$, $(II)$, $(III)$, and $(IV)$ are \dq essentially equivalent'' (equivalent modulo those multiplicative factors). The implication $(I)\Rightarrow(IV)$ is probably the most usesful of the lot. If $f$ satisfies $(IV)$ and has integral 0, we can decompose it as a rapidly converging series of H\"older continuous, compactly supported functions having cancelation, with good control on the functions' sizes and supports \cite[Lemma 3.5]{U}\cite[Lemma 6.3]{Wilsonbook}. Such decompositions have value in their own right.


\medskip

\begin{thebibliography}{99}

\bibitem{AbramSilv} F. Abramovich, B. W. Silverman, \dq Wavelet Decomposition Approaches to Statistical Inverse Problems'', \emph{Biometrika} {\bf 85} (1998), 115--129.\medskip

\bibitem{CaiPlanning} T. T. Cai, \dq On adaptive wavelet estimation of a derivative and other related linear inverse problems'', \emph{Journal of Statistical Planning and Inference} {\bf 108} (2002), 329--349.\medskip

\bibitem{WavCoifMeyer} R. R. Coifman, Y. Meyer, \emph{Wavelets: \Ca-Zygmund and multilinear operators}, Cambridge University Press (1997), Cambridge.\medskip

\bibitem{DidJafPip} G. Didier, S. Jaffard, V. Pipiras, \dq On the vaguelet and Riesz properties of $L^2$-unbounded transformations of orthogonal wavelet bases'', \emph{Journal of Approximation Theory} {\bf 176} (2013), 94--117.\medskip

\bibitem{Jia}R.-Q. Jia, \dq Bessel sequences in Sobolev spaces'', \emph{Applied and Computational Harmonic Analysis} {\bf 20} (2006), 298-311.\medskip

\bibitem{LiuYY}H. Liu, H. Yang, Q. Yang, \dq Carleson Measures and Trace Theorem for $\beta$-harmonic Functions'', \emph{Taiwanese Journal of Mathematics} {\bf 22} (2018), 1107--1138.\medskip

\bibitem{U} A. Uchiyama, \dq A constructive proof of the Fefferman-Stein decomposition of $BMO(\R^n)$'', \emph{Acta Mathematica} {\bf 148} (1982), 215--241.\medskip

\bibitem{Wilsonbook}M. Wilson, \emph{Weighted Littlewood-Paley Theory
and Exponential-Square Integrability}, Springer Lecture Notes in
Mathematics {\bf 1924} (2007), New York.

\end{thebibliography}
\end{document}